\documentclass[a4paper,10pt,reqno]{amsart}
\usepackage{marginnote,amsmath,amsfonts,amscd,amssymb,graphicx,mathrsfs,eufrak}
\usepackage[dvips,all,arc,curve,color,frame]{xy}
\usepackage[usenames]{color}

\usepackage[colorlinks]{hyperref}

\newtheorem{thm}{Theorem}[section]
\newtheorem{prop}[thm]{Proposition}
\newtheorem{lemma}[thm]{Lemma}
\newtheorem{defin}[thm]{Definition}
\newtheorem{cor}[thm]{Corollary}

\theoremstyle{definition}

\newtheorem{remark}[thm]{Remark}
\newtheorem{conj}[thm]{Conjecture}

\newcommand{\C}[1]{\mathbb{C}^{#1}}
\newcommand{\m}{\mathfrak{m}}

\newcommand{\p}{\mathfrak{p}}

\renewcommand{\t}[1]{\textnormal{#1}}

\newcommand{\looptop}[2]{\xy \SelectTips{cm}{10}
\POS(0,0) \endxy}

\newcommand{\blob}{{\scriptscriptstyle\bullet}}

\def\op{\mathop{\rm op}\nolimits}

\def\CM{\mathop{\rm CM}\nolimits}

\def\depth{\mathop{\rm depth}\nolimits}

\def\hgt{\mathop{\rm ht}\nolimits}
\def\mod{\mathop{\rm mod}\nolimits}

\def\refl{\mathop{\rm ref}\nolimits}
\def\proj{\mathop{\rm proj}\nolimits}

\def\pd{\mathop{\rm proj.dim}\nolimits}

\def\Hom{\mathop{\rm Hom}\nolimits}
\def\End{\mathop{\rm End}\nolimits}
\def\Ext{\mathop{\rm Ext}\nolimits}

\def\add{\mathop{\rm add}\nolimits}
\def\Cok{\mathop{\rm Cok}\nolimits}

\def\ker{\mathop{\rm ker}\nolimits}

\def\sup{\mathop{\rm sup}\nolimits}

\def\Spec{\mathop{\rm Spec}\nolimits}
\def\Max{\mathop{\rm Max}\nolimits}
\def\gl{\mathop{\rm gl.dim}\nolimits}

\CompileMatrices
\MakeOutlines

\edef\marginnotetextwidth{\the\textwidth}

\newcommand{\new}[1]{{\color{blue} #1}}

\begin{document}
\title[\textsc{On the NC Bondal--Orlov Conjecture}]{\textsc{On the Noncommutative Bondal--Orlov Conjecture}}
\author{Osamu Iyama}
\address{Osamu Iyama\\ Graduate School of Mathematics\\ Nagoya University\\ Chikusa-ku, Nagoya, 464-8602, Japan.}
\email{iyama@math.nagoya-u.ac.jp}
\author{Michael Wemyss}
\address{Michael Wemyss, School of Mathematics, James Clerk Maxwell Building, The King's Buildings, Mayfield Road, Edinburgh, EH9 3JZ, UK.}
\email{wemyss.m@googlemail.com}
\begin{abstract}
Let $R$ be a normal, equi-codimensional Cohen-Macaulay ring of dimension $d\geq 2$ with a canonical module $\omega_R$.   We give a sufficient criterion that establishes a derived equivalence between the noncommutative crepant resolutions of $R$.  When $d\leq 3$ this criterion is always satisfied and so all noncommutative crepant resolutions of $R$ are derived equivalent.
Our method is based on cluster tilting theory for commutative algebras, developed in \cite{IW}.
\end{abstract}
\maketitle
\parindent 20pt
\parskip 0pt

\section{Introduction}

The following conjecture, due to Bondal--Orlov, is one of the main motivating problems in the study of derived categories in higher-dimensional birational geometry.

\begin{conj}[{Bondal--Orlov}] If
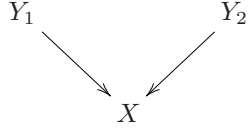
\[
\xymatrix{
Y_1\ar[rd]&&Y_2\ar[ld]\\
&X
}
\]
are two crepant resolutions of $X$, then $Y_1$ is derived equivalent to $Y_2$.
\end{conj}

In the study of one-dimensional fibres \cite{VdB1d}, and also in the McKay Correspondence for dimension $d\leq 3$ \cite{BKR}, $Y_1$ and $Y_2$ are derived equivalent to certain noncommutative rings.  Thus in these cases, showing that $Y_1$ and $Y_2$ are derived equivalent is equivalent to showing that the corresponding noncommutative rings are derived equivalent.   With this motivation, these noncommutative structures were axiomized by Van den Bergh \cite{VdBNCCR} into the concept of a noncommutative crepant resolution (=NCCR):

\begin{defin}\label{nonsingorder}
Let $R$ be a Cohen-Macaulay (=CM) ring and let $\Lambda$ be a module finite $R$-algebra.\\
\t{(1)} {\rm \cite{Aus78, Auslanderisolated, CR90}} $\Lambda$ is called an \emph{$R$-order} if $\Lambda$ is a maximal Cohen-Macaulay $R$-module.
An $R$-order $\Lambda$ is called \emph{non-singular} if $\gl \Lambda_\p=\dim R_\p$ for all primes $\p$ of $R$.\\
\t{(3)}  {\rm \cite{VdBNCCR}} By a \emph{noncommutative crepant resolution} of $R$ we mean $\Gamma:=\End_R(M)$ where $M$ is a non-zero reflexive $R$-module such that $\Gamma$ is a non-singular $R$-order.
\end{defin}

We remark here that Van den Bergh defined NCCRs only when the base ring $R$ is Gorenstein, since these are the types of varieties which have a chance of admitting crepant resolutions and so there is a good analogy with the geometry.  However when $R$ is CM and non-Gorenstein there are sometimes many NCCRs of $R$, and these are related to cluster tilting objects in the category $\CM R$ \cite[5.5]{IW}.  Thus, although geometrically we are only really interested in NCCRs when $R$ is Gorenstein, there are strong algebraic reasons to consider the more general case.

\begin{conj}[{\cite[4.6]{VdBNCCR}}]\label{VdBconj} If $R$ is a normal Gorenstein domain, then all crepant resolutions of $\Spec R$ (commutative and noncommutative) are derived equivalent.
\end{conj}
We remark that, even in dimension three, it is still not known in full generality whether the existence of a commutative crepant resolution is equivalent to the existence of a NCCR.  Certainly in dimension four and higher this fails in both directions; there are examples where there exist NCCRs but no commutative crepant resolutions (e.g.\ any terminal cyclic quotient singularity in dimension four) and also examples where there exist crepant resolutions but no NCCRs (e.g.\ $R=\C{}[[x_0,x_1,x_2,x_3,x_4]]/(x_0^5+x_1^4+x_2^4+x_3^4+x_4^4)$, see \cite[3.5]{Dao}).

Nevertheless there are still no known counterexamples to \ref{VdBconj}.  Because of the breakdown in the correspondence between geometry and algebra outlined in the above paragraph, in this paper we restrict ourselves to the following special case of \ref{VdBconj}:

\begin{conj}[{Noncommutative Bondal--Orlov}]\label{NCBOconj} If $R$ is a normal Gorenstein domain, then all NCCRs of $R$ are derived equivalent.
\end{conj}

Again we remark that $R$ is assumed to be Gorenstein only so that there is an analogy with the geometry; if we are brave enough then we can hope that \ref{NCBOconj} is true in the more general setting of CM rings. Surprisingly, we show below (\ref{main1intro} and \ref{mainintro}) that this hope is not unfounded.

In recent work \cite{IW} we argued that, at least in dimension three, the condition of finite global dimension in the definition of NCCR should be weakened and so we defined and studied the notion of a maximal modification algebra (MMA).  Conjecturally this corresponds to studying  maximal crepant partial resolutions of $\Spec R$ which may or may not be smooth.  The work \cite{IW} is quite technical since we had to overcome problems involving infinite global dimension, but the tools we used and developed were quite general and not necessarily restricted to the setting considered in \emph{loc.\ cit}.

The main purpose of this paper is to illustrate that when we work in the full generality of $R$ being CM, but restrict ourselves to studying NCCRs (instead of MMAs), the techniques in \cite{IW} give very simple, short proofs of some relatively strong results.  The first surprising result is that in dimension three the Noncommutative Bondal--Orlov Conjecture holds in the full generality of CM rings.  This generalizes \cite[8.8]{IR} to algebras which do not have Gorenstein base rings:

\begin{thm}\label{main1intro}(=\ref{d=2},\ref{main1})
Let $R$ be a $d$-dimensional CM equi-codimensional (see \ref{equicodim}) normal domain with a canonical module $\omega_R$. \\
\t{(1)} If $d=2$, then all NCCRs of $R$ are Morita equivalent.\\
\t{(2)} If $d=3$, then all NCCRs of $R$ are derived equivalent.
\end{thm}

In fact the above follows immediately from a more general result.  Below in \ref{depthleveld} and \ref{mainintro} we give a sufficient condition for arbitrary dimension $\dim R=d$ to establish when any two given NCCRs of $R$ are derived equivalent.  This condition is empty when $d\leq 3$, and so \ref{mainintro} proves \ref{main1intro}.  Unlike other conditions in the literature (e.g.\ the sufficient condition in \cite{BKR}) our condition is checked on the base singularity $R$, not on a fibre product, and is thus readily verified either by hand or by using computer algebra.

\begin{defin}\label{depthleveld}
Let $R$ be a normal, equi-codimensional CM ring of dimension $d\geq 3$, with a canonical module $\omega_{R}$. 
Suppose $\Lambda:=\End_{R}(M)$ and $\End_R(N)$ are NCCRs of $R$.  Consider a complex (see \ref{construction of complex})
\begin{eqnarray}
0\to M_{d-2}\to \cdots\to M_{1}\to M_{0}\to N\to0\label{downintro}
\end{eqnarray}
with each $M_{i}\in\add M$ such that applying $\Hom_R(M,-)$ induces a projective resolution
\[
0\to \Hom_R(M,M_{d-2})\to \cdots\to \Hom_R(M,M_{0})\to \Hom_R(M,N)\to 0
\]
of $\Hom_R(M,N)$ as a $\Lambda$-module.  We say that the pair $(M,N)$ satisfies the \emph{depth condition} if 
\[
\depth_{R_{\m}}\Hom_{R}(M_{i},N)_{\m}\geq d-i-1
\]
for all $\m\in\Max R$ and all $i\ge0$. Note that this inequality is automatic for $i=d-3$ and $i=d-2$ (see \ref{depthofhom}).
\end{defin}

Note that to check whether $(M,N)$ satisfies the depth condition, one only needs to know the modules $M_0,\hdots,M_{d-4}$, not the whole complex.  For unexplained terminology and more explanation, see \S\ref{prelim} and \S\ref{mainsection}.  We also remark here that the complex (\ref{downintro}) is built by taking successive right $\add M$-approximations on the base singularity $R$ and so requires very little explicit knowledge of the rings $\End_R(M)$ or $\End_R(N)$.  See \ref{addapproxremark} for more details.

Our main result is the following, where $(-)^*=\Hom_R(-,R)$.

\begin{thm}\label{mainintro}(=\ref{higherdimNCBO})
Let $R$ be a normal, equi-codimensional CM ring of dimension $d\geq 2$, with canonical module $\omega_{R}$.  Suppose that $\End_{R}(M)$ and $\End_{R}(N)$ are NCCRs such that both $(M,N)$ and $(N^{*},M^{*})$ satisfy the depth condition. Then:\\
\t{(1)} $\Hom_{R}(M,N)$ is a tilting $\End_{R}(M)$-module of projective dimension at most $d-2$.\\
\t{(2)} $\End_{R}(M)$ and $\End_{R}(N)$ are derived equivalent.
\end{thm}

\medskip
\noindent
{\bf Conventions.}
All modules are left modules, so for a ring $A$ we denote $\mod A$ to be the category of finitely generated left $A$-modules.  Throughout when composing maps $fg$ will mean $f$ then $g$.  Note that with this convention $\Hom_R(M,X)$ is a $\End_R(M)$-module and $\Hom_R(X,M)$ is a $\End_R(M)^{\rm op}$-module.  For $M\in\mod A$ we denote $\add M$ to be the full subcategory consisting of summands of finite direct sums of copies of $M$, and we denote $\proj A:=\add A$ to be the category of finitely generated projective $A$-modules.  Throughout we will always use the letter $R$ to denote some kind of \emph{commutative noetherian} ring. 

\medskip
\noindent
{\bf Acknowledgement.}
The authors thank Ryo Takahashi for useful information on the Acyclicity Lemma.

\section{Preliminaries}\label{prelim}

In this section we record the commutative algebra preliminaries that are required in the proof of the main theorem, and also fix notation.  The main technical results we need are \ref{reflequiv} and \ref{ABlocal}, both of which already appear in the literature.

For a commutative noetherian local ring $(R,\m)$ and $M\in\mod R$, recall that the \emph{depth} of $M$ is defined to be
\[
\depth_R M:=\inf \{ i\geq 0: \Ext^i_R(R/\m,M)\neq 0 \},
\]
which coincides with the maximal length of a $M$-regular sequence.  Keeping the assumption that $(R,\m)$ is local, we say that $M\in\mod R$ is \emph{maximal Cohen-Macaulay} (or simply, \emph{CM}) if $\depth_R M=\dim R$.  This definition generalizes to the non-local case as follows: if $R$ is an arbitrary commutative noetherian ring, we say that $M\in\mod R$ is \emph{CM} if $M_\p$ is CM for all prime ideals $\p$ in $R$. We say that $R$ is a \emph{CM ring} if $R$ is a CM $R$-module.  We denote $\CM R$ to be the category of CM $R$-modules.

The following is well-known and is easy to check from the definition.
\begin{lemma}[{see e.g.\ \cite[2.3]{IW}}]\label{depthofhom}
Let $(R,\m)$ be a local ring of dimension $d\geq 2$ and let $\Lambda$ be a module finite $R$-algebra. Then for all $M,N\in\mod \Lambda$ with $\depth_{R} N\geq 2$, we have $\depth_{R}\Hom_{\Lambda}(M,N)\geq 2$.
\end{lemma}

Throughout this paper we denote
\[(-)^{*}:=\Hom_{R}(-,R):\mod R\to\mod R\]
and we say that $X\in\mod R$ is \emph{reflexive} if the natural map $X\to X^{**}$ is an isomorphism.  We denote $\refl R$ to be the category of reflexive $R$-modules. When $R$ is a normal domain, the category $\refl R$ is closed under both kernels and extensions.  Note that if $\depth R\geq 2$ and $X\in\refl R$, then $\depth_R X\geq 2$ by applying \ref{depthofhom} to $X=\Hom_R(X^*,R)$.  

In this paper (see \ref{reflequiv}) reflexive modules over noncommutative rings play a crucial role.  

\begin{defin}\label{reflexivedefinition} Let $R$ be any commutative ring. If $\Lambda$ is any $R$-algebra then we say that $M\in\mod \Lambda$ is a \emph{reflexive $\Lambda$-module} if it is reflexive as an $R$-module.  
\end{defin}
Note that we do not require that the natural map $M\to \Hom_{\Lambda^{\op}}(\Hom_\Lambda(M,\Lambda),\Lambda)$ is an isomorphism. When the underlying commutative ring $R$ is a normal domain, the following reflexive equivalence is crucial: 
\begin{lemma}\label{reflequiv}
Suppose $R$ is a commutative ring and $M\in\mod R$. \\
\t{(1)} The functor $\Hom_{R}(M,-):\mod R\to\mod\End_{R}(M)$ restricts to an equivalence 
\[ 
\Hom_{R}(M,-):\add M\stackrel{\approx}{\to}\proj\End_{R}(M).
\]
\t{(2)} If $R$ is a normal domain and $M\in\refl R$ is non-zero, then we have an equivalence
\[
\Hom_{R}(M,-):\refl R\stackrel{\approx}{\to}\refl\End_{R}(M).
\]
\end{lemma}
\begin{proof}
(1) is famous `projectivization' (e.g. \cite[II.2.1]{ARS}).\\
(2) is also well-known (e.g. \cite[2.4]{IR}, \new{\cite{RV89}}).
\end{proof}

For some global--local arguments to work we require the following weak assumption, which we note is automatically satisfied in the main examples we are interested in, namely affine domains \cite[13.4]{Eisenbud}.
\begin{defin}\label{equicodim}
A commutative ring $R$ is called \emph{equi-codimensional} if all its maximal ideals have the same height.
\end{defin}

Now let $R$ be an equi-codimensional CM ring of dimension $d$ with a canonical module $\omega_{R}$.  Recall that for a non-local CM ring $R$, a finitely generated $R$-module $\omega_{R}$ is called a \emph{canonical module} if $({\omega_{R}})_{\m}$ is a canonical $R_{\m}$-module for all $\m\in\Max R$ \cite[3.3.16]{BH}.  In this case $(\omega_{R})_{\p}$ is a canonical $R_{\p}$-module for all $\p\in\Spec R$ since canonical modules localize for local CM rings \cite[3.3.5]{BH}.

When $R$ is local, we have the following Auslander--Buchsbaum type equailty, which in particular says that the $\Lambda$-modules which are CM as $R$-modules are precisely the projective $\Lambda$-modules.

\begin{lemma}\label{ABlocal}
Let $(R,\m)$ be a local CM ring with a canonical module $\omega_{R}$, and let $\Lambda$ be a non-singular $R$-order.  Then for any $X\in\mod \Lambda$,
\[
\depth_R X+\pd_{\Lambda} X=\dim R.
\]
\end{lemma}
\begin{proof}
\cite[2.3]{IR}, or combine \cite[2.13, 2.14]{IW}.
\end{proof}

As an application, let us prove the following observation.

\begin{prop}\label{d=2}
Let $R$ be a $2$-dimensional CM equi-codimensional normal domain with a canonical module $\omega_R$.  If $R$ has a NCCR, then all NCCRs of $R$ are Morita equivalent.
\end{prop}
\begin{proof}
Let $\End_{R}(M)$ be a NCCR.
By \ref{reflequiv}, we have an equivalence
\[\Hom_R(M,-):\refl R\to\refl\End_R(M).\]
Any $X\in\refl\End_R(M)$ satisfies $\depth_{R_\m}X_\m\geq2$ for all $\m\in\Max R$ by \ref{depthofhom}, hence $X_\m\in\proj\End_{R_\m}(M_\m)$ for all $\m\in\Max R$ by \ref{ABlocal}. Thus $X$ is a projective $\End_R(M)$-module, and so $\refl\End_R(M)=\proj\End_R(M)$. By \ref{reflequiv}(1), this implies that $\refl R=\add M$.

If $\End_R(N)$ is another NCCR, then by above we have $\add M=\refl R=\add N$.  Thus $\End_{R}(M)$ and $\End_{R}(N)$ are Morita equivalent, via the progenerator $\Hom_R(M,N)$.
\end{proof}

One of the benefits of the existence of a canonical module is that depth localizes well:
\begin{lemma}\label{depthlocalizesCAN}
Let $(R,\m)$ be a local CM ring with canonical module $\omega_R$, and let $X\in\mod R$.  
For all $\p\in\Spec R$ we have
\[
\dim R-\depth_RX\geq\dim R_\p-\depth_{R_\p}X_\p.
\]
\end{lemma}
\begin{proof}
It is well known (see e.g.\ \cite[3.5.11]{BH}) that we have equalities 
\begin{eqnarray*}
\dim R-\depth_{R}X &=& \sup\{ i\geq 0 \mid \Ext^i_{R}(X,\omega_R)\neq 0\}\\
\dim R_\p-\depth_{R_\p}X_\p &=& \sup\{ i\geq 0 \mid \Ext^i_{R_\p}(X_\p,\omega_{R\p})\neq 0\}.
\end{eqnarray*}
Since $\Ext^i_{R_\p}(X_\p,\omega_{R\p})=\Ext^i_{R}(X,\omega_R)_\p$, we have the assertion.
\end{proof}

Finally, for the convenience of the reader we recall the definition of Miyashita's tilting module:

\begin{defin}\label{gentiltingdef}\cite{H88,M86}
Let $\Lambda$ be a \emph{noetherian} ring.  Then $T\in\mod \Lambda$ is called a \emph{tilting module} if the following conditions are satisfied:
\begin{itemize}
\item $\pd_{\Lambda}T<\infty$,
\item $\Ext^{i}_{\Lambda}(T,T)=0$ for all $i>0$,
\item There exists an exact sequence
\[
0\to \Lambda\to T_{0}\to T_{1}\to\cdots \to T_{t-1}\to T_{t}\to 0
\]
with each $T_{i}\in\add T$.
\end{itemize}
\end{defin}

It is a classical result in tilting theory \cite{H88} that $\Lambda$ and $\End_\Lambda(T)$ are derived equivalent, for any tilting $\Lambda$-module $T$.

\section{Proof of Main Theorem}\label{mainsection}

Throughout this section we let $R$ be a $d$-dimensional CM equi-codimensional normal domain with $d\geq 2$ and canonical module $\omega_R$. 
Suppose $\Lambda:=\End_{R}(M)$ is a NCCR.

\begin{lemma}\label{construction of complex}
For any $N\in\refl R$, there exists a complex 
\begin{eqnarray}
0\to M_{d-2}\to \cdots\to M_{1}\to M_{0}\to N\to0\label{approxres1}
\end{eqnarray}
with each $M_{i}\in\add M$ such that applying $\Hom_R(M,-)$ induces a projective resolution
\begin{eqnarray}
0\to \Hom_R(M,M_{d-2})\to \cdots\to \Hom_R(M,M_{0})\to \Hom_R(M,N)\to 0\label{projres1}
\end{eqnarray}
of $\Hom_R(M,N)$ as a $\Lambda$-module.
\end{lemma}
\begin{proof}
We have $\depth_{R_{\m}}\Hom_R(M,N)_{\m}\geq 2$ for all $\m\in\Max R$ by \ref{depthofhom}, and so by \ref{ABlocal} $\pd_{\Lambda_{\m}}\Hom_R(M,N)_{\m}\leq d-2$ for all $\m\in\Max R$.  Hence $\pd_{\Lambda} \Hom_R(M,N)\leq d-2$.  Consequently, we have the assertion by \ref{reflequiv}.
\end{proof}

\begin{defin}
With the assumptions as above, in particular $\dim R=d$, then we say that the pair $(M,N)$ satisfies the \emph{depth condition} if there exists a complex (\ref{approxres1}) such that
\[
\depth_{R_{\m}}\Hom_{R}(M_{i},N)_{\m}\geq d-i-1
\]
for all $\m\in\Max R$ and all $i\geq0$.  
Note that this inequality is automatic for $i=d-3$ and $i=d-2$ (see \ref{depthofhom}).
\end{defin}

When $d\leq 3$ the condition is empty and so every such $(M,N)$ satisfies the depth condition.  

\begin{remark}\label{addapproxremark}
To construct a complex (\ref{approxres1}) involves very little explicit knowledge of $\End_R(M)$ or $\End_R(N)$, since it is built by taking successive right $\add M$-approximations on the base ring $R$.  
Recall that for $R$-modules $M$ and $N$, we say that
a morphism $f:M_0\to N$ is a \emph{right $\add M$-approximation}
if $M_0\in \add M$ and further the map
\[
\Hom_R(M,M_0)\stackrel{\cdot f}{\to}\Hom_R(M,N)
\] 
is surjective. This just says that $f$ is a morphism into $N$ from a module in $\add M$ such that any other morphism into $N$ from any other module in $\add M$ factors through $f$. 
\end{remark}

We require the following technical result.

\begin{lemma}\label{vanishing}
Let $R$ be a local CM ring of dimension $d$, let
\[
0\stackrel{f^{-1}}{\to} X^0\xrightarrow{f^0}X^1\xrightarrow{f^1}\cdots\xrightarrow{f^{d-2}}X^{d-1}\xrightarrow{f^{d-1}}X^d\xrightarrow{f^d}X^{d+1}\xrightarrow{f^{d+1}}\cdots
\]
be a complex of finitely generated $R$-modules, and let $H^i$ be the $i$-th homology.  Assume that $\depth X^i\ge d-i$ for all $i\ge0$, and $H^i=0$ for all $i\ge d$.\\
\t{(1)} If $H^i$ is a finite length $R$-module for all $i\ge0$, then the above sequence is exact.\\
\t{(2)} If $H^i_\p=0$ for all $i\ge0$ and all $\p\in\Spec R$ with $\hgt\p\le i$, then the above sequence is exact.
\end{lemma}

\begin{proof}
The point is that any finite length module with positive depth has to be zero.\\
(1)  Although this is known as the Acyclicity Lemma (e.g. \cite[1.4.23]{BH}), we give a proof for the convenience of the reader.  We show that every $H^j$ is zero by induction on $j$. 
Since $0\to H^0=\ker f^0\to X^0$ is exact, by the depth lemma $\depth H^0>0$ and so since $H^0$ has finite length, necessarily $H^0=0$.  Thus the result is true for $j=0$, so fix $0\le j<d$ and assume that $H^k=0$ for all $0\le k<j$.
Then we have exact sequences
\[
0\to X^0\xrightarrow{f^0}\cdots\xrightarrow{f^{j-1}}X^j\to\Cok f^{j-1}\to0,
\]
\[
0\to H^j\to\Cok f^{j-1}\to X^{j+1}.
\]
From the first sequence we have $\depth\Cok f^{j-1}\ge d-j>0$ by using depth lemma repeatedly.
Since $H^j$ is a finite length $R$-module, the second sequence forces $H^j=0$.\\
(2) Let $\p$ be minimal among prime ideals of $R$ satisfying $H^i_\p\neq0$ for some $i$.
We have a complex
\begin{equation}\label{localized sequence}
0\to X^0_\p\xrightarrow{f^0_\p}X^1_\p\xrightarrow{f^1_\p}\cdots\xrightarrow{f^{d-2}_\p}X^{d-1}_\p\xrightarrow{f^{d-1}_\p}X^d_\p\xrightarrow{f^d_\p}X^{d+1}_\p\xrightarrow{f^{d+1}_\p}\cdots
\end{equation}
of finitely generated $R_\p$-modules.
By \ref{depthlocalizesCAN}, we have
\[
\depth_{R_\p}X^i_\p\ge\dim R_\p+\depth_RX^i-\dim R\ge \hgt\p+(d-i)-d=\hgt\p-i
\]
for all $i\geq 0$.  We have $H^i_\p=0$ for all $i\ge\hgt\p$ by our assumption. Moreover $H^i_\p$ is a finite length $R_\p$-module for all $i\ge0$ by our choice of $\p$.
Thus \eqref{localized sequence} is exact by (1), a contradiction.
Hence $H^i=0$ for all $i\ge0$.
\end{proof}

Now we are ready to prove our main result.

\begin{thm}\label{higherdimNCBO}
Let $R$ be a normal, equi-codimensional CM ring of dimension $d\geq 2$, with a canonical module $\omega_{R}$.  Suppose that $\End_{R}(M)$ and $\End_{R}(N)$ are NCCRs such that both $(M,N)$ and $(N^{*},M^{*})$ satisfy the depth condition.  Then:\\
\t{(1)} $\Hom_{R}(M,N)$ is a tilting $\End_{R}(M)$-module of projective dimension at most $d-2$.\\
\t{(2)} $\End_{R}(M)$ and $\End_{R}(N)$ are derived equivalent.
\end{thm}
\begin{proof}
Set $\Lambda:=\End_{R}(M)$ and $T:=\Hom_{R}(M,N)$. We show that $T$ is a tilting $\Lambda$-module with $\pd_{\Lambda}T\leq d-2$.  The fact that $\End_{\Lambda}(T)\cong\End_{R}(N)$ follows from \ref{reflequiv}(2).\\
(i) We have $\pd_{\Lambda}T\leq d-2$ by \ref{construction of complex}. \\
(ii)  We now show that $\Ext_\Lambda^i(T,T)=0$ for all $i>0$.  Let us recall the notation:
\[
0\to M_{d-2}\to\cdots\to M_1\to M_0\to N\to0
\]
is a complex which induces a projective resolution
\[
0\to\Hom_R(M,M_{d-2})\to\cdots\to\Hom_R(M,M_0)\to\Hom_R(M,N)\to0.
\]
Applying $\Hom_\Lambda(-,T)$ and using reflexive equivalence \ref{reflequiv}(2) we obtain a complex
\[
X^\blob:=\quad0\to\End_R(N)\to\Hom_R(M_0,N)\to\cdots\to\Hom_R(M_{d-2},N)\to0
\]
that satisfies the following properties:
\begin{itemize}
\item Let $X^0:=\End_R(N)$ and $X^i:=\Hom_R(M_{i-1},N)$ for $i\ge1$. Then for all $\m\in\Max R$, we have $\depth_{R_\m} X^0_\m=d$  since $X^0\in\CM R$, and  $\depth_{R_\m} X^i_\m\ge d-i$ for all $i\ge1$ by the depth condition.
\item Let $H^i$ be the homology of the complex $X^\blob$ at $X^i$.  Then $H^0=H^1=0$ and $H^i=\Ext^{i-1}_\Lambda(T,T)$ for all $i\ge2$.
\item Since $\pd_{\Lambda_\p}T_\p\le\hgt\p-2$ by \ref{depthofhom} and \ref{ABlocal}, we have $H^i_\p=\Ext^{i-1}_{\Lambda_\p}(T_\p,T_\p)=0$ for all $i\ge0$ and all $\p\in\Spec R$ with $\hgt\p\le i$.
\end{itemize}
Thus for every $\m\in\Max R$ if we consider the complex $X^\blob_\m$ of $R_\m$-modules, all the assumptions in \ref{vanishing}(2) are satisfied, so the sequence $X^\blob_\m$ is exact.
Thus we have $\Ext^i_\Lambda(T,T)_\m=0$ for all $i>0$ and all $\m\in\Max R$.  This implies that $\Ext^i_\Lambda(T,T)=0$ for all $i>0$.\\
(iii) Applying the same argument as in (i) and (ii) to $(N^*,M^*)$, we have a complex
\[0\to N_{d-2}^*\to\cdots\to N_1^*\to N_0^*\to M^*\to0\]
inducing an exact sequence
\[
0\to\End_R(M^*)\to\Hom_R(N_0^*,M^*)\to\cdots\to\Hom_R(N_{d-2}^*,M^*)\to 0.
\]
Since $(-)^{*}:\refl R\to \refl R$ is a duality, this means that
\[
0\to\End_{R}(M)\to\Hom_{R}(M,N_{0})\to\cdots\to\Hom_{R}(M,N_{d-2})\to 0
\]
is exact.   But this is simply 
\[
0\to\Lambda\to T_{0}\to T_{1}\to\cdots\to T_{d-2}\to 0
\]
with all $T_{i}\in \add T$.  Hence $T$ is a tilting $\Lambda$-module with $\pd_{\Lambda}T\leq d-2$.
\end{proof}

Since the depth condition is automatic for $d\leq 3$, the next result follows immediately from \ref{higherdimNCBO}.

\begin{cor}\label{main1}
Let $R$ be a $3$-dimensional CM equi-codimensional normal domain with a canonical module $\omega_R$.  If $R$ has a NCCR, then all NCCRs of $R$ are derived equivalent.
\end{cor}

\end{document}